\documentclass[letterpaper,11pt]{article}
\usepackage{stmaryrd}
\usepackage{cite}
\usepackage[all]{xy}
\usepackage[dvips]{graphicx}
\usepackage{graphpap}
\usepackage{ifthen}
\usepackage{enumerate}
\usepackage{amssymb}
\usepackage{mathrsfs}
\usepackage[errorshow]{tracefnt}
\usepackage{fancybox}
\usepackage{amsfonts}
\usepackage{amsmath}
\usepackage{amssymb}
\usepackage{multirow}
\usepackage{array}
\usepackage{mathtools}
\usepackage{soul}
\usepackage{pict2e}
\usepackage{mathabx}
\DeclareMathAlphabet{\mathpzc}{OT1}{pzc}{m}{it}

\makeatletter
\newcommand{\thickhline}{%
    \noalign {\ifnum 0=`}\fi \hrule height 1pt
    \futurelet \reserved@a \@xhline
}
\newcolumntype{'}{@{\hskip\tabcolsep\vrule width 1pt\hskip\tabcolsep}}
\makeatother
\newcolumntype{"}{@{\hskip\tabcolsep\vrule width 1.5pt\hskip\tabcolsep}}
\makeatother

\def\boxit#1{\vbox{\hrule\hbox{\vrule\kern3pt
             \vbox{\kern3pt#1\kern3pt}\kern3pt\vrule}\hrule}}

\newcommand{\beq}{\begin{equation}}
\newcommand{\beqn}{\begin{equation*}}
\newcommand{\eeq}{\end{equation}}
\newcommand{\eeqn}{\end{equation*}}
\newcommand{\beqa}{\begin{eqnarray}}
\newcommand{\beqan}{\begin{eqnarray*}}
\newcommand{\eeqa}{\end{eqnarray}}
\newcommand{\eeqan}{\end{eqnarray*}}
\newcommand{\bdm}{\begin{displaymath}}
\newcommand{\edm}{\end{displaymath}}

\newcommand{\alg}{\mathcal G}
\newcommand{\ssub}{\mathcal H} %\mathfrak{h}}
\newcommand{\scom}{\mathcal E} %E}%{\mathfrak{p}}
\newcommand{\sub}{{\scriptscriptstyle{\mathcal H}}} %\mathfrak{h}}
\newcommand{\com}{{\scriptscriptstyle{\mathcal E}}}%{\mathfrak{p}}

\newcommand{\ba}{\begin{array}}
\newcommand{\ea}{\end{array}}

\newcommand\nn{\nonumber}

\newcommand\benu{\begin{enumerate}}
\newcommand\eenu{\end{enumerate}}
\newcommand\bit{\begin{itemize}}
\newcommand\eit{\end{itemize}}

\newtheorem{theorem}{Theorem}[section]

\newtheorem{cor}[theorem]{Corollary}
\newtheorem{prop}[theorem]{Proposition}

\def\Pf{\noindent \textbf{Proof. }}
\def\End{\mathrm{End\,}}
\def\Hom{\mathrm{Hom\,}}

\def\der'{\mathfrak{der}'\,}
\def\der{\mathfrak{der}\,}
\def\str'{\mathfrak{str}'\,}
\def\str{\mathfrak{str}\,}

\def\so{\mathfrak{so}}

\def\qed{\hspace{\stretch{1}} $\Box$ \\
\noindent}

\newcommand{\dlb}{\llbracket}%{\ensuremath{[\![}}
\newcommand{\drb}{\rrbracket}%{\ensuremath{]\!]}}

\newcommand{\blb}%{\ensuremath
{\text{$\llbracket$\hspace{-4pt}\scalebox{0.99}{$|$}\hspace{-2.58pt}\scalebox{0.99}{$|$}\hspace{-2.58pt}\scalebox{0.99}{$|$}}\,\!}
\newcommand{\brb}%{\ensuremath{
{\text{$\rrbracket$\hspace{-4pt}\scalebox{0.99}{$|$}\hspace{-2.58pt}\scalebox{0.99}{$|$}\hspace{-2.58pt}\scalebox{0.99}{$|$}}}

\def\*{\partial}

\RequirePackage{hyperref} 

\numberwithin{equation}{section}

\begin{document}

\frenchspacing

\includegraphics[height=1.85cm]{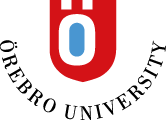}

\vskip-59.5pt
\hfill {\tt 2012.10954v3}\\ %\today} \\
\vskip-10pt
\hfill {\tt\today}\\ %\tt \jobname} \\

\vspace*{4cm}

\begin{center}
\noindent
{\LARGE {\sf \textbf{Nonlinear realisations of Lie superalgebras}}}\\
\vspace{.3cm}

\renewcommand{\thefootnote}{\fnsymbol{footnote}}

\vskip 1truecm

\noindent
{\large {\sf \textbf{Jakob Palmkvist}
}}\\
\vskip .5truecm
        {\it 
        {
        School of Science and Technology\\\"Orebro University\\ SE-701 82 \"Orebro, Sweden}\\[3mm]}
        {\tt jakob.palmkvist@oru.se} \\
\end{center}

\vskip 1cm

\centerline{\sf \textbf{
Abstract}}
\vskip .2cm

\noindent
For any decomposition of a Lie superalgebra $\alg$ into a direct sum $\alg=\ssub\oplus \scom$ 
of a subalgebra $\ssub$ and a subspace $\scom$, without any further restrictions on $\ssub$ and $\scom$,
we construct a nonlinear realisation of
$\alg$ on $\scom$. The result generalises a theorem by Kantor from Lie algebras to Lie superalgebras.
When $\alg$ is a differential graded Lie algebra,
we show that it gives a construction of an associated $L_\infty$-algebra.
%\underline{\it Keywords:} Lie superalgebras, Leibniz algebras, differential graded Lie algebras, $L_\infty$ algebras.

\noindent

\newpage

\pagestyle{plain}

\tableofcontents

\section{Introduction}

Representations 
of Lie algebras can be generalised to nonlinear realisations.
This means that 
the elements in the Lie algebra are mapped to operators that are not necessarily linear, but constant,
quadratic or of higher order. In many important applications, the operators act on a vector space which can be identified with a subspace of
the Lie algebra itself, complementary to a subalgebra.
One example is the conformal realisation of the Lie algebra $\so(2,D)$ on a $D$-dimensional vector space,
based on the decomposition of $\so(2,D)$ as a 3-graded Lie algebra $G=G_{-1}\oplus G_0\oplus G_1$,
where $G_0=\so(1,D-1)$ and $G_{\pm1}$ are $D$-dimensional subspaces. % such that $[G_1,G_{-1}]\subseteq G_0$.
In this realisation, the subalgebra $\so(1,D-1)$ acts linearly, whereas the two $D$-dimensional subspaces can be considered as consisting of constant
and quadratic operators, respectively.
In other examples, the linearly realised subalgebra is not the degree-zero subalgebra in a $\mathbb{Z}$-grading, but 
defined by being pointwise invariant under an involution.

For any
decomposition of a Lie algebra $G$ into a direct sum $G=H \oplus E$ of a subalgebra $H$ and a complementary subspace $E$,
there is formula for a nonlinear realisation of $G$ on $E$ given by Kantor \cite{Kantor}.
The conformal realisation of 
a semisimple Lie algebra with a 3-grading $G=G_{-1}\oplus G_0 \oplus G_1$
is obtained 
from this formula in the special case where 
$H=G_0 \oplus G_1$ and $E = G_{-1}$.
The corresponding application to a
semisimple
Lie algebra
with a 5-grading $G=G_{-2}\oplus G_{-1}\oplus G_0 \oplus G_1 \oplus G_2$ leads to a
{\it quasiconformal} realisation if the subspaces $G_{\pm2}$ are one-dimensional
\cite{Gunaydin:2000xr,Palmkvist:2005gc}.
In these cases the subspace $E$ is actually also a subalgebra, but this need not be the case in the general formula.
There are no restrictions on $[E,E]$, nor on $[H,E]$; the only requirement is $[H,H]\subseteq H$.

In this paper, we generalise Kantor's formula from Lie algebras to Lie superalgebras.
Also in the restriction to Lie algebras, our proof is very different from Kantor's, being
purely algebraic, without references to homogeneous spaces for Lie groups. 

We expect our generalisation to be useful 
in applications to physics, in particular to models where a Lie superalgebra
can be used to organise
the field content or to encode the gauge structure.
In such cases it might be interesting to investigate whether the Lie superalgebra also can be realised as a symmetry.
We also expect the result to be relevant for applications 
of other related structures, such as Leibniz algebras, differential graded Lie algebras and $L_\infty$-algebras, for which a renewed interest has appeared recently 
in the context of gauge theories, see for example Refs.\
\cite{Lavau:2017tvi,Kotov:2018vcz,Bonezzi:2019bek,Lavau:2019oja,Bonezzi:2019ygf,Lavau:2020pwa,Borsten:2021ljb}.
In fact, our framework illuminates the relations between these structures. 
In particular,
our main result leads to the construction of an $L_\infty$-algebra
associated to any differential graded Lie algebra.

The paper is organised as follows.
\begin{itemize}
\item We start in Section 2 with an arbitrary vector space $U_1$. 
We associate a $\mathbb{Z}$-graded Lie algebra $U$ to it, from which we in turn construct the
Lie algebra $S$ of symmetric operators on $U_1$. 
The $\mathbb{Z}$-graded Lie algebra $U$ associated to a vector space
$U_1$ was introduced in Ref.~\cite{Kantor-graded}, but here we use a different recursive approach, following Refs.~\cite{Palmkvist:2005gc,Palmkvist:2009qq}.

\item In Section 3, we modify the construction: we then start with a vector space $\mathcal U_1$ that is equipped with a $\mathbb{Z}_2$-grading,
to which we associate a $\mathbb{Z}$-graded Lie superalgebra $\mathcal{U}$ \cite{Palmkvist:2009qq}. From $\mathcal{U}$ we
construct the Lie superalgebra $\mathcal{S}$ of symmetric operators on $\mathcal{U}_1$ (where the symmetry is now
actually a $\mathbb{Z}_2$-graded symmetry).

\item In Section 4, we furthermore assume that the vector space $\mathcal U_1$ itself
is a Lie superalgebra $\mathcal G$.
This means that it
is equipped with a Lie superbracket, consistent with the
$\mathbb{Z}_2$-grading already present in Section 3.
We show that it extends to a Lie superbracket on $\mathcal S$ (different from the one defined in Section 3).

\item In Section 5, we still assume $\mathcal U_1=\alg$
but also that this Lie superalgebra decomposes into a direct sum $\alg = \ssub \oplus \scom$,
where $\ssub$ is a subalgebra. We show that it extends to a corresponding direct sum $\mathcal S = \mathcal S_\sub \oplus \mathcal S_\com$.
As our main result, Theorem 5.4, we show that there is a Lie superalgebra homomorphism from $\mathcal G$ to
$\mathcal S_\com$. This result generalises the main theorem in Ref.\ \cite{Kantor} from Lie algebras to Lie superalgebras.

\item In Section 6, we assume that $\mathcal U_1=\alg$ itself has a $\mathbb{Z}$-grading consistent with the $\mathbb{Z}_2$-grading,
and is equipped with a differential, turning it into a differential graded Lie algebra.
As an example of an application of our main result, we use it in order to
construct an $L_\infty$-algebra from $\alg$, and show that the brackets agree with those given explictly in Ref.\ \cite{Getzler:2010}.

\end{itemize}

\noindent
\underline{\it Acknowledgments:}
I would like to thank Martin Cederwall, Sylvain Lavau and Arne Meurman for discussions. I am particularly grateful to Sylvain Lavau, who have also given many useful comments
on the manuscript.

\section{The \texorpdfstring{$\mathbb{Z}$}{Z}-graded Lie algebra associated to a vector space}

We start with an arbitrary vector space $U_1$ over some field of characteristic zero, from which
we define vector spaces $U_0,U_{-1},U_{-2},\ldots$ recursively by
\begin{align} \label{defUsubspaces}
U_{-p+1}=\Hom(U_1,U_{-p+2})
\end{align} 
for $p=1,2,\ldots$. Thus $U_{-p+1}$ consists of all linear maps from $U_1$ to $U_{-p+2}$, and in particular $U_0 = \End{U_1}$.

Let $A_p \in U_{-p+1}$, for some $p=1,2,\ldots$, and let $x_1,x_2,\ldots
\in U_1$. Then $A_p(x_1) \in U_{-p+2}$ and if $p\geq 2$,
this means that $A_p(x_1)(x_2)=\big(A_p(x_1)\big)(x_2)$
is an element in $U_{-p+3}$. Continuing in this way, we finally find that
$A(x_1)(x_2)\cdots(x_p)$
is an element in $U_1$, which we may also write as $A(x_1,x_2,\ldots,x_p)$.
Thus we have a vector space isomorphism
\begin{align}
U_{-p+1}=\Hom(U_1,U_{-p+2}) \simeq \Hom((U_1)^p,U_1)
\end{align}
and we may consider elements in $U_{-p+1}$ not only as linear maps from $U_1$ to $U_{-p+2}$
but also as linear maps from $(U_1)^p$ to $U_1$, or as $p$-linear operators on $U_1$.  
We will refer to elements in $U_{-p+1}$ simply as operators of order $p$, even for $p=0$, so that
the elements in $U_1$ are considered as operators of order zero.

\subsection{The Lie algebra \texorpdfstring{$U_{0-}$}{U0-}}

Next we let $U_{0-}$ be the direct sum of the vector spaces defined in the previous section,
$U_{0-}=U_0 \oplus U_{-1} \oplus U_{-2} \oplus \cdots$.
For any
$A_p \in U_{-p+1}$ (where $p=1,2,\ldots$)
and any $x \in U_1$, we write
\begin{align} \label{ring-initial}
A_p \circ x &= A_p(x)\,, & x \circ A_p &= 0\,.
\end{align}
We then define a map
\begin{align}
\circ \quad :\quad U_{-p+1} \times U_{-q+1} \to U_{-(p+q-1)+1}\quad ,\quad (A_p,B_q) \mapsto A_p \circ B_q
\end{align}
for any $p,q =1,2,\ldots$ recursively by
\begin{align} \label{circrec}
(A \circ B)(x) = A \circ B(x) + A(x) \circ B
\end{align}
and extend it to a bilinear operation on $U_{0-}$ by linearity.
For $p=q=1$ this is the usual composition of (linear) maps,
\begin{align} \label{linmapcomp}
(A_1 \circ B_1)(x)=A_1 \circ B_1(x) + A_1(x) \circ B_1 = A_1 \big(B_1(x)\big)\,,
\end{align}
where the last equality follows from (\ref{ring-initial}) since $B_1(x)$ and $A_1(x)$ are elements in $U_1$.
We give two more examples,
\begin{align}
(A_2 \circ B_1)(x_1,x_2)&=(A_2 \circ B_1)(x_1)(x_2)\nn\\
&=\big( A_2 \circ B_1(x_1) + A_2(x_1) \circ B_1 \big)(x_2)\nn\\
&=\big(A_2 \circ B_1(x_1)\big) (x_2) + A_2(x_1)\circ B_1(x_2)+A_2(x_1)(x_2) \circ B_1 \nn\\
&=A_2 \big(  B_1(x_1)\big) (x_2) + A_2(x_1)\big(B_1(x_2)\big)\nn\\
&=A_2 \big(  B_1(x_1),x_2\big) + A_2\big(x_1,B_1(x_2)\big)\,,
\end{align}
\begin{align}
(B_1 \circ A_2)(x_1,x_2)&=(B_1 \circ A_2)(x_1)(x_2)\nn\\
&=\big( B_1 \circ A_2(x_1) + B_1(x_1) \circ A_2 \big)(x_2)\nn\\
&=\big( B_1 \circ A_2(x_1) \big)(x_2)\nn\\
&=B_1 \circ A_2(x_1)(x_2)+B_1(x_2) \circ A_2(x_1)\nn\\
&=B_1 \circ A_2(x_1,x_2)=B_1\big(A_2(x_1,x_2)\big)\,,
\end{align}
which 
are
easily 
generalised to 
\begin{align}
(A_p \circ B_1)(x_1,x_2,\ldots,x_p)&=A_p\big(B_1(x_1),x_2,\ldots,x_p)\big)\nn\\
&\quad\,+A_p\big(x_1,B_1(x_2),\ldots,x_p)\big)\nn\\
&\quad\,+\cdots+A_p\big(x_1,x_2,\ldots,B_1(x_p))\big)\,,
\end{align}
\begin{align}
(B_1 \circ A_p)(x_1,x_2,\ldots,x_p)=B_1\big(A_p(x_1,x_2,\ldots,x_p)\big)\,.
\end{align}
In these examples, the subscripts of the operators indicate their orders (whereas the subscripts on elements $x$ in $U_1$ are
just labels used to distinguish them from each other).

The property $ U_{-p+1} \circ U_{-q+1} \subseteq U_{-(p+q-1)+1}$ means that $U_{0-}$ is a $\mathbb{Z}$-graded algebra with respect to $\circ$
(with vanishing subspaces corresponding to positive integers). The following proposition
says that this algebra furthermore is associative.

\begin{prop} \label{assalg}
The vector space $U_{0-}$ together with the bilinear operation $\circ$ is an associative algebra.
\end{prop}
\Pf We will show that
\begin{align} \label{associativity}
\big((A_p \circ B_q) \circ C_r\big)(x) = \big(A_p \circ (B_q \circ C_r)\big)(x)
\end{align}
for any triple of operators $A_p,B_q,C_r$ of order $p,q,r\geq1$, respectively, and any $x \in U_1$. We do this by induction over $p+q+r\geq 3$.
When $p+q+r=3$, we have $p=q=r=1$, and the assertion follows by (\ref{linmapcomp}). Suppose now that it holds when $p+q+r=s$ for some $s\geq3$,
and set $p+q+r=s+1$. We then have (omitting the subscripts)
\begin{align}
\big(A \circ (B \circ C)\big)(x)&=A \circ (B\circ C)(x) + A(x) \circ (B \circ C)\nn\\
&=A \circ \big(B\circ C(x)\big)+ A \circ \big(B(x) \circ C\big) + A(x) \circ (B \circ C)\nn\\
&=(A \circ B)\circ C(x)+ \big(A \circ B(x)\big) \circ C + \big(A(x) \circ B\big) \circ C\nn\\
&=(A \circ B)\circ C(x)+ (A \circ B)(x) \circ C \nn\\
&=\big((A \circ B)\circ C\big)(x)\,,
\end{align}
using the induction hypothesis in the third step,
and the proposition follows by the principle of induction.
\qed

\noindent
Note that the identity (\ref{associativity}) is not satisfied when $r=0$ and $p,q \neq0$. Then (omitting the subscripts and setting $C_0=x$) we instead have the (right) Leibniz identity
\begin{align}
(A \circ B) \circ x = A \circ (B \circ x)+ (A \circ x)\circ B\,.
\end{align}

For any
$A,B \in U_{0-}$, we now set
\begin{align} \label{dbracketdef}
\llbracket A,B \rrbracket = A \circ B - B \circ A
\end{align}
and we have the following obvious consequence of Proposition \ref{assalg}.
\begin{cor} \label{cor:u-lie}
The vector space $U_{0-}$ together with the bracket $\llbracket \cdot,\cdot\rrbracket$ is a
Lie algebra.
\end{cor}

\subsection{Extending \texorpdfstring{$U_{0-}$}{U0-} to \texorpdfstring{$U$}{U}}

Let $U_+=U_1 \oplus U_2 \oplus \cdots$ be the free Lie algebra generated by the vector space $U_1$ (with the natural $\mathbb{Z}_+$-grading)
and set
\begin{align}
U=U_{0-} \oplus  U_+ = \cdots \oplus  U_{-1} \oplus  U_0 \oplus  U_1 \oplus
 U_2 \oplus \cdots\,.
\end{align}
We will use the notation
\begin{align}
 U_{i+}&=\bigoplus_{j\geq i} U_j\,, &  U_{i-}&=\bigoplus_{j\leq i} U_j\,
\end{align}
for any $i\in \mathbb{Z}$. 

We use the same notation $\llbracket \cdot,\cdot\rrbracket$ for the
two Lie brackets on $U_{0-}$ and $U_+$, respectively, and we will now unify them into one Lie bracket on the whole of $U$, the direct sum of these two vector spaces.
We thus have to define brackets $\llbracket A,u \rrbracket =- \llbracket u,A \rrbracket$ for any $A \in U_{0-}$ and $u \in U_+$.
For $u=x \in U_1$ we set 
\begin{align} \label{bracketmedx}
\llbracket A,x \rrbracket=A(x)\,. 
\end{align}
If $u \in U_{2+}$, then we may assume that
$u=\llbracket v,w \rrbracket$ for some $v,w \in U_+$.
We then define recursively
\begin{align}
\llbracket A,\llbracket v,w\rrbracket \rrbracket=\llbracket \llbracket A,v\rrbracket,w\rrbracket-\llbracket \llbracket A,w\rrbracket,v\rrbracket\,\label{defAvw}
\end{align}
and extend the bracket by linearity to the case when $u$ is a sum of such terms $\llbracket v,w \rrbracket$.
In order to ensure that the definition is meaningful, we have to show that it respects the Jacobi identity
in the sense that
\begin{align}
\llbracket A,\llbracket \llbracket u,v\rrbracket,w\rrbracket \rrbracket =\llbracket A,\llbracket u,\llbracket v,w\rrbracket \rrbracket \rrbracket 
-\llbracket A,\llbracket v,\llbracket u,w\rrbracket \rrbracket \rrbracket \,
\end{align}
for any $A \in U_{0-}$ and $u,v,w\in U_+$. Indeed, we get
\begin{align}
\llbracket A,\llbracket \llbracket u,v\rrbracket ,w\rrbracket \rrbracket 
&=\llbracket \llbracket A,\llbracket u,v\rrbracket \rrbracket ,w\rrbracket -\llbracket \llbracket A,w\rrbracket ,\llbracket u,v\rrbracket \rrbracket \nn\\
&=\llbracket \llbracket \llbracket A,u\rrbracket ,v\rrbracket ,w\rrbracket -\llbracket \llbracket \llbracket A,v\rrbracket ,u\rrbracket ,w\rrbracket 
-\llbracket \llbracket \llbracket A,w\rrbracket ,u\rrbracket ,v\rrbracket \rrbracket +\llbracket \llbracket \llbracket A,w\rrbracket ,v\rrbracket ,u\rrbracket \rrbracket \nn\\
&=\llbracket \llbracket \llbracket A,u\rrbracket ,w\rrbracket ,v\rrbracket -\llbracket \llbracket \llbracket A,v\rrbracket ,w\rrbracket ,u\rrbracket 
-\llbracket \llbracket \llbracket A,w\rrbracket ,u\rrbracket ,v\rrbracket \rrbracket +\llbracket \llbracket \llbracket A,w\rrbracket ,v\rrbracket ,u\rrbracket \rrbracket \nn\\
&\quad\,+\llbracket \llbracket A,u\rrbracket ,\llbracket v,w\rrbracket \rrbracket -\llbracket \llbracket A,v\rrbracket ,\llbracket u,w\rrbracket \rrbracket \nn\\
&=\llbracket \llbracket A,\llbracket u,w\rrbracket \rrbracket ,v\rrbracket -\llbracket \llbracket A,\llbracket v,w\rrbracket \rrbracket ,u\rrbracket \nn\\
&\quad\,+\llbracket \llbracket A,u,\rrbracket\llbracket v,w\rrbracket \rrbracket -\llbracket \llbracket A,v\rrbracket ,\llbracket u,w\rrbracket \rrbracket \nn\\
&=\llbracket A,\llbracket u,\llbracket v,w\rrbracket \rrbracket \rrbracket -\llbracket A,\llbracket v,\llbracket u,w\rrbracket \rrbracket \rrbracket \,
\end{align}
using Jacobi identities like
\begin{align}
\llbracket \llbracket A,w\rrbracket ,\llbracket u,v\rrbracket \rrbracket =\llbracket \llbracket \llbracket A,w\rrbracket ,u\rrbracket ,v\rrbracket \rrbracket -\llbracket \llbracket \llbracket A,w\rrbracket ,v\rrbracket ,u\rrbracket \rrbracket \,.
\end{align}
Such Jacobi identities
follow either (if $\llbracket A,w\rrbracket\in U_+$)
by the fact that $U_+$ is a Lie algebra or (if $\llbracket A,w\rrbracket\in U_{0-}$) by the definition (\ref{defAvw}).
\begin{prop}
The vector space $U=U_{0-} \oplus U_+$ together with the bracket $\llbracket \cdot,\cdot\rrbracket$ is a Lie algebra.
\end{prop}
\Pf
The Jacobi identities with either all three elements in $U_{0-}$ or all three elements in $U_{+}$ are satisfied, by Corollary \ref{cor:u-lie} and by the
construction of $U_+$ as a free Lie algebra. Also the Jacobi identities with one element in $U_{0-}$ and two elements in $U_+$ are satisfied,
by the definition (\ref{defAvw}). It only remains to check the Jacobi identities with two elements $A,B \in U_{0-}$ and one element $u\in U_+$.
Assuming that $u$ is homogeneous with respect to the $\mathbb{Z}$-grading, $u \in U_k$,
this can be done by induction over $k\geq 1$. For $k=1$, we have
\begin{align}
\dlb \dlb A,B \drb,u\drb &= \dlb A,B \drb(u)=(A \circ B - B\circ A)(u)\nn\\
&=A\circ B(u)+A(u)\circ B - B\circ A(u)-B(u)\circ A\nn\\
&=\dlb A, B(u)\drb-\dlb B, A(u)\drb=\dlb A,\dlb B,u\drb\drb-\dlb B,\dlb A,u\drb\drb\,.
\end{align}
For $k\geq 2$, we may (as above), assume that $u=\dlb v,w\drb$, where $v,w \in U_{1+}$. Assuming furthermore (as the induction hypothesis) that
\begin{align}
\dlb \dlb A,B \drb,v\drb &= \dlb A,\dlb B,v\drb\drb-\dlb B,\dlb A,v\drb\drb\,,\nn\\
\dlb \dlb A,B \drb,w\drb &= \dlb A,\dlb B,w\drb\drb-\dlb B,\dlb A,w\drb\drb\,,
\end{align}
we get
\begin{align}
\llbracket \llbracket A,B\rrbracket ,u\rrbracket  &= \llbracket \llbracket A,B\rrbracket ,\llbracket v,w\rrbracket \rrbracket \nn\\
&=\llbracket \llbracket \llbracket A,B\rrbracket ,v\rrbracket ,w\rrbracket -\llbracket \llbracket \llbracket A,B\rrbracket ,w\rrbracket ,v\rrbracket \nn\\
&=\llbracket \llbracket A,\llbracket B,v\rrbracket \rrbracket ,w\rrbracket -\llbracket \llbracket B,\llbracket A,v\rrbracket \rrbracket ,w\rrbracket 
-\llbracket \llbracket A,\llbracket B,w\rrbracket \rrbracket ,v\rrbracket 
+\llbracket \llbracket B,\llbracket A,w\rrbracket \rrbracket ,v\rrbracket \nn\\
&=\llbracket A,\llbracket \llbracket B,v\rrbracket ,w\rrbracket \rrbracket -\llbracket B,\llbracket \llbracket A,v\rrbracket ,w\rrbracket \rrbracket 
-\llbracket A,\llbracket \llbracket B,w\rrbracket ,v\rrbracket \rrbracket +\llbracket \llbracket B,\llbracket \llbracket A,w\rrbracket ,v\rrbracket \rrbracket \nn\\
&\quad\,-\llbracket \llbracket B,v\rrbracket ,\llbracket A,w\rrbracket \rrbracket +\llbracket \llbracket A,v\rrbracket ,\llbracket B,w\rrbracket \rrbracket 
+\llbracket \llbracket B,w\rrbracket ,\llbracket A,v\rrbracket \rrbracket -\llbracket \llbracket A,w\rrbracket ,\llbracket B,v\rrbracket \rrbracket \nn\\
&=\llbracket A,\llbracket B,\llbracket v,w\rrbracket \rrbracket \rrbracket -\llbracket B,\llbracket A,\llbracket v,w\rrbracket \rrbracket \rrbracket \nn\\
&=\llbracket A,\llbracket B,u\rrbracket \rrbracket -\llbracket B,\llbracket A,u\rrbracket \rrbracket \,,
\end{align}
and it follows by induction that also these Jacobi identities are satisfied.
\qed

\subsection{The Lie algebra \texorpdfstring{$S$}{S} of symmetric operators}

For $p\geq 0$, let $S_{-p+1}$ be the subspace of $U_{-p+1}$ consisting of 
elements $A_p \in U_{-p+1}$ such that
$\llbracket A_p,u \rrbracket \subseteq U_{2+}$ for all $u \in U_{2+}$, and set
\begin{align}
S=\bigoplus_{p\geq0} S_{-p+1}\,.
\end{align}
Because of the $\mathbb{Z}$-grading, if $p\geq2$ then $S_{-p+1}$ consists of all operators $A_p$ such that
$\llbracket A,u \rrbracket=0$ for all $u \in U_{2+}$, whereas $S_{0}=U_0$ and $S_1=U_1$.
Furthermore, $S \oplus U_{2+}$ is the idealiser (or normaliser) of $U_{2+}$ in $U$, and $S$ can be identified with
the quotient space obtained by factoring out $U_{2+}$ from its idealiser in $U$.

We will refer to elements in $S_{-p+1}$ as {\it symmetric} operators of order $p$, and a linear combination of
symmetric operators will also be called a symmetric operator, even if it is not homogenous with respect to the
$\mathbb{Z}$-grading.
Note that we consider all elements in $U_0$ as symmetric operators of order
one, and even all elements in $U_1$ as symmetric operators of order zero.

It follows easily by the Jacobi identity that if $A$ is a symmetric operator of order one or higher, then $A \circ x$ is a symmetric operator
as well, for any $x \in U_1$.

The operators in $U$ of order two or higher included in $S$ are indeed precisely those that are symmetric in the following sense.
If $A_2 \in S_{-1}$ and $x,y\in U_1$ (so that $\llbracket x,y\rrbracket \in U_2$), then
\begin{align}
0=\llbracket A_2, \llbracket x,y\rrbracket\rrbracket
&=\llbracket \llbracket  A_2, x\rrbracket,y\rrbracket-
\llbracket \llbracket  A_2, y\rrbracket,x\rrbracket \nn\\
&=A_2(x)(y)-A_2(y)(x)=A_2(x,y)-A_2(y,x)\,,
\end{align}
so that $A_2(x,y)=A_2(y,x)$. It is straightforward to show that generally, the condition
$\llbracket A_p,u \rrbracket=0$ for all $u \in U_{2+}$ is equivalent to the condition that 
\begin{align}
A_p(x_1,\ldots,x_p)=A_p(y_1,\ldots,y_p)\,,
\end{align} 
where $(y_1,\ldots,y_p)$ is any permutation of $(x_1,\ldots,x_p)$. We write this (as usual) as
\begin{align} \label{symmetrinotation}
A_p(x_1,\ldots,x_p)=A_p(x_{(1},\ldots,x_{p)})\,,
\end{align} 
where the right hand side denotes $1/p!$ times the sum of $A(y_1,\ldots,y_p)$ over all permutations $(y_1,\ldots,y_p)$ of $(x_{1},\ldots,x_{p})$.

For any symmetric operator $A_p$ there is a unique corresponding map $U_1 \to U_1$ (non-linear if $p\neq1$) given by
$x \mapsto A(x,x,\ldots,x)$.
In order to characterise a symmetric operator $A_p$ it is thus sufficient to set $x_1=x_2=\cdots =x_p$ in
$A(x_1,\ldots,x_p)$.

In particular for symmetric operators, it is convenient to 
replace the 
bilinear operation $\circ$ on $U_{0-}$
by another one, which 
differs from $\circ$ by normalisation.
We define
a bilinear operation 
$\bullet$ on $U_{0-}$ by
\begin{align} \label{bulletdef}
A_p \bullet B_q = \frac{p!\,q!}{(p+q-1)!}A_p \circ B_q\,
\end{align}
for $A_p \in U_{-p+1}$ and $B_q \in U_{-q+1}$.
If $A_p$ and $B_q$ are symmetric operators, we then get
\begin{align}
(A_p \bullet B_q)
(x,x,\ldots,x)=
pA_p(B_q(x,x,\ldots,x),x,\ldots,x)\,.
\end{align}
Since $A_p \circ B_q$ is an operator of order $p+q-1$, the 
linear map $\phi : U_{0-} \to U_{0-}$ given by
\begin{align}
\phi(A_p) = \frac1{p!}A_p
\end{align}
satisfies
\begin{align}
\phi(A) \bullet \phi(B) = \phi(A \circ B)
\end{align}
for any two operators $A$ and $B$
and thus the two algebras obtained by equipping the vector space $U_{0-}$ with
$\circ$ and $\bullet$, respectively, are isomorphic to each other.

We now extend the bilinear operation $\bullet$ from $U_{0-}$ to $U_{1-}$. First we set
\begin{align}
A_p \bullet x &= p (A_p \circ x) = pA_p(x)\,, & x \bullet A_p &=0  
\end{align}
for $A_p \in U_{-p+1}$ (where $p=1,2,\ldots)$ and $x\in U_1$. Thus the definition
(\ref{bulletdef}) is still valid if we allow one of $A_p$ and $B_q$ to be an operator of order zero,
that is, an element in $U_1$.
For example, we have
\begin{align}
(A_p \bullet y)(x_1,\ldots, x_{p-1})=p A_p(y,x_1,\ldots, x_{p-1})\,,
\end{align}
whereas
\begin{align}
(A_p \circ y)(x_1,\ldots, x_{p-1})= A_p(y,x_1,\ldots, x_{p-1})\,.
\end{align}
Second, we set
\begin{align}
x \bullet y =0
\end{align}
for $x,y \in U_1$ in order to close $U_{1-}$ under $\bullet$. This makes the operation $\bullet$ really different from $\circ$
(not only up to normalisation),
since we kept $x \circ y$ undefined.

For any $A,B \in S$, we set
\begin{align}
\blb A,B \brb = A \bullet B - B \bullet A\,.
\end{align}
It follows that the vector space $S$ equipped with this bracket is a Lie algebra isomorphic to
the quotient algebra obtained by factoring out $U_{2+}$ from the idealiser of $U_{2+}$ in $U$. Moreover, if
$U_1$ is $n$-dimensional, it is straightforward to show
that $S$ is isomorphic to the Lie algebra $W_n$ of formal vector fields $\sum_{i=1}^n f_i \tfrac{\partial}{\partial x_i}$, where
$f_i$ are formal power series in $n$ variables $x_1,\ldots, x_n$.  

\section{Generalisation from Lie algebras to Lie superalgebras}

We will now 
repeat the steps in the preceding section in a more general case. Instead of starting with an arbitrary vector space $U_1$ we now start with
an arbitrary $\mathbb{Z}_2$-graded vector space $\mathcal U_1$.
Thus $\mathcal U_1$ can be
decomposed into a direct sum $\mathcal U_1=\mathcal U_{1}{}^{(0)} \oplus \mathcal U_{1}{}^{(1)}$ of 
two subspaces $\mathcal U_{1}{}^{(0)}$ and $\mathcal U_{1}{}^{(1)}$. Like for any $\mathbb{Z}_2$-graded vector space, these subspaces
(and their elements) are said to be {\it even} and {\it odd},
respectively.
This leads to a corresponding decomposition
\begin{align}
\mathcal U_{p+1}=\mathcal U_{p+1}{}^{(0)}\oplus \mathcal U_{p+1}{}^{(1)}
\end{align}
of each vector space $\mathcal U_{p+1}$, by refining (\ref{defUsubspaces}) to
\begin{align}
\mathcal U_{p+1}{}^{(0)}&=\Hom(\mathcal U_1{}^{(0)},\mathcal U_{-p+2}{}^{(0)})\oplus \Hom(\mathcal U_1{}^{(1)},\mathcal U_{-p+2}{}^{(1)}) \,,\nn\\
\mathcal U_{p+1}{}^{(1)}&=\Hom(\mathcal U_1{}^{(0)},\mathcal U_{-p+2}{}^{(1)})\oplus \Hom(\mathcal U_1{}^{(1)},\mathcal U_{-p+2}{}^{(0)}) \,.
\end{align}
Now, let $\mathcal U_+=\mathcal U_1 \oplus \mathcal U_2 \oplus \cdots$ be the free {Lie superalgebra} generated by $\mathcal U_1$ (with the natural $\mathbb{Z}$-grading) and set
\begin{align}
\mathcal U=\mathcal U_{0-} \oplus \mathcal U_+ = \cdots \oplus \mathcal U_{-1} \oplus \mathcal U_0 \oplus \mathcal U_1 \oplus \mathcal U_2 \oplus \cdots\,.
\end{align}
We thus have a $\mathbb{Z}_2$-graded vector space
$\mathcal U = \mathcal U^{(0)} \oplus \mathcal U^{(1)}$.
If $u \in \mathcal U_{0-}{}^{(i)}$ for $i=0,1$, we 
use the notation $|u|=i$ for the
$\mathbb{Z}_2$-degree of $u$.

We can now repeat the steps in the preceding section,
carrying over notation and terminology in a straightforward way. The formulas will however
differ from those in the preceding section by factors of powers of $(-1)$, where we (without loss of generality) have to assume that the elements in $\mathcal U$ that appear are
homogeneous with respect to the $\mathbb{Z}_2$-grading.

Thus, we equip $\mathcal U_{0-}$ with an associative bilinear operation $\circ$, from which we define a Lie superbracket $\dlb\cdot,\cdot\drb$ on 
$\mathcal U_{0-}$. In these definitions, we modify 
(\ref{circrec}) to
\begin{align}
(A \circ B)(x) = A \circ B(x) + (-1)^{|B||x|} A(x) \circ B
\end{align}
and (\ref{dbracketdef}) to
\begin{align}
\llbracket A,B \rrbracket = A \circ B - (-1)^{|A||B|} B \circ A\,.
\end{align}
When we then unify the
brackets on $\mathcal{U}_{0-}$ and $\mathcal{U}_+$ to one on the whole of $\mathcal U$,
we keep the definition 
$\llbracket A,x \rrbracket=A(x)$ in (\ref{bracketmedx}),
but modify (\ref{defAvw}) to
\begin{align}
\llbracket A,\llbracket v,w\rrbracket \rrbracket=\llbracket \llbracket A,v\rrbracket,w\rrbracket
-(-1)^{|v||w|}\llbracket \llbracket A,w\rrbracket,v\rrbracket\,.
\end{align}

We do not give the proofs here, since they differ from those given in the preceding section only by
factors of
powers of $(-1)$.

In fact, the modifications made here are actually generalisations, since
the Lie superalgebra $\mathcal U$
reduces to the original Lie algebra $U$ in the special case where $\mathcal U_1$ has a trivial odd subspace, $\mathcal U_1{}^{(1)}=0$.
Thus, starting with a vector space $\mathcal U_1$, we can decompose it in different ways into a direct sum of an even and an odd subspace, which lead to different associated $\mathbb{Z}$-graded Lie superalgebras $\mathcal U$.
The decomposition where $\mathcal U_1$ is considered as an even vector space (coinciding with its even subspace) leads to the associated $\mathbb{Z}$-graded Lie algebra described in the preceding section.
But we can also consider it as an odd vector space. Only in this case the $\mathbb{Z}$-grading of the Lie superalgebra
is {\it consistent}, which means that $\mathcal U_i \subseteq \mathcal U_{(0)}$ if $i$ is even and $\mathcal U_i \subseteq \mathcal 
U_{(1)}$ if $i$ is odd.
 
The Lie superalgebra $\mathcal S$, constructed from $\mathcal U$ in the same way as $S$ is constructed from $U$, now consists of operators with a $\mathbb{Z}_2$-graded symmetry, rather than purely symmetric ones. However, for simplicity we will still refer to them as symmetric operators.
Generalising the notation (\ref{symmetrinotation}), we denote $\mathbb{Z}_2$-graded symmetry with angle brackets rather than ordinary parentheses, so that
\begin{align}
A_p(x_{1},\ldots,x_{p})=A_p(x_{\langle 1},\ldots,x_{p\rangle})
\end{align}
if $A_p \in \mathcal S$, 
where the right hand side denotes $1/p!$ times the sum of $(-1)^\varepsilon A(y_1,\ldots,y_p)$ over all permutations $(y_1,\ldots,y_p)$ of
$(x_{1},\ldots,x_{p})$, where $\varepsilon$ is the number of
transpositions of two odd elements.
 
\subsection{Leibniz algebras}

In the next section we will assume that the $\mathbb{Z}_2$-graded vector space $\mathcal U_1$ is a Lie superalgebra. Before that,
we will briefly give another example of a case where $\mathcal U_1$ is an algebra. In many such cases, identities for elements in this 
algebra can be reformulated
as identities for elements in the associated $\mathbb{Z}$-graded Lie superalgebra $\mathcal U$,
including the bilinear operation of the algebra as an element in $\mathcal U_{-1}$.

A (left) {\it Leibniz algebra} is an algebra $\mathcal U_1$ where the bilinear operation $\odot$ satisfies the (left) Leibniz identity
\begin{align}
x \odot (y \odot z) = (x \odot y) \odot z + y \odot (x \odot z)\,.
\end{align}
If we now consider $\mathcal U_1$ as a $\mathbb{Z}_2$-graded vector space with trivial even subspace
and let $\Theta$ be the element in $\mathcal U_{-1}$ associated to $\odot$ by 
\begin{align} \label{leibniz}
x \odot y = \Theta(x,y) = \Theta(x)(y) = \dlb\dlb \Theta ,x \drb,y\drb
\end{align}
then the Leibniz identity (\ref{leibniz}) is equivalent to the condition 
\begin{align}
\dlb \Theta,\Theta \drb=0\,. 
\end{align}
(Since $\mathcal U_1$ is odd, $\mathcal U_{-1}$ is odd as well, and the condition $\dlb \Theta,\Theta \drb=0$ is not trivially satisfied,
but equivalent to $\Theta \circ \Theta =0$). 
Indeed, 
by the Jacobi identity (keeping in mind that $\Theta,x,y,z$ are all odd),
\begin{align}
\dlb \Theta,\Theta \drb(x,y,z) &= \dlb \dlb \dlb \dlb \Theta,\Theta \drb,x\drb,y\drb,z\drb\nn\\
&=2\dlb \dlb \dlb \Theta,\dlb \Theta ,x \drb\drb,y\drb,z\drb\nn\\
&=2\dlb \dlb  \Theta,\dlb\dlb \Theta ,x \drb,y\drb\drb,z\drb+2 \dlb \dlb  \dlb\Theta,y\drb,\dlb \Theta ,x \drb\drb,z\drb\nn\\
&=2\dlb \dlb  \Theta,\dlb\dlb \Theta ,x \drb,y\drb\drb,z\drb+2 \dlb \dlb  \Theta,y\drb,\dlb \dlb\Theta ,x \drb,z\drb\drb-2 \dlb \dlb  \Theta,x\drb,\dlb \dlb\Theta ,y \drb,z\drb\drb\nn\\
&=2(x \odot y) \odot z +2 y \odot (x \odot z) -2 x \odot (y \odot z) \,.
\end{align}
Now let $\langle \Theta \rangle$ be the one-dimensional subspace of $\mathcal U_{-1}$ spanned by $\Theta$. Since $\dlb \Theta,\Theta \drb=0$,
the subspace $\langle \Theta \rangle \,\oplus\, \mathcal U_{0+}$ of $\mathcal U$ is a subalgebra. This Lie superalgebra can also be considered as 
a differential graded Lie algebra $\mathcal U_{0+}$ with a differential $\dlb \Theta,\cdot \drb$. Thus any Leibniz algebra gives rise to a differential graded Lie algebra \cite{Lavau:2017tvi,Kotov:2018vcz,Lavau:2019oja}. In Section 6 we will see how in turn any
differential graded Lie algebra gives rise to an $L_\infty$-algebra.

\section{The case when \texorpdfstring{$\mathcal U_1$}{U1} is a Lie superalgebra \texorpdfstring{$\alg$}{G}}

We now assume not only that $\mathcal U_1$ is a $\mathbb{Z}_2$-graded vector space, but furthermore that 
$\mathcal U_1$ is a Lie superalgebra $\alg$ 
with a bracket $[\cdot,\cdot]$. We 
extend the bracket to 
the whole of $\mathcal U_{1-}$ recursively by
\begin{align} \label{bracketdef}
[A,B] \bullet x = [A,B\bullet x]+(-1)^{|x||B|}[A \bullet x,B]\,.
\end{align}
We recall that any operator $A_p$ of order $p$ is defined by its action on $\mathcal U_1$, and that $A_p \bullet x = p A_p(x)$.
If $B=y$ is an element in $\mathcal U_0$, that is, an operator of order zero, then $y \bullet x=0$, so that
\begin{align}
[A,y]\bullet x = [A,y\bullet x]+(-1)^{|x||y|}[A\bullet x,y]=(-1)^{|x||y|}[A\bullet x,y]\,.
\end{align}
\begin{prop} \label{u1-isaliealgebra}
The $\mathbb{Z}_2$-graded vector space $\mathcal U_{1-}$ together with the bracket $[\cdot,\cdot]$ is a Lie superalgebra.
\end{prop}
\Pf
We will show that the Jacobi identity
\begin{align}
[[A,B],C] = [A,[B,C]]+(-1)^{|B||C|}[[A,C],B]
\end{align}
is satisfied for any triple of operators $A,B,C$ of order $p,q,r$, respectively, by induction over $p+q+r\geq0$. When $p+q+r=0$, we have $p=q=r=0$
and the Jacobi identity is satisfied since $\alg$ is a Lie superalgebra. If we assume that it is satisfied when $p+q+r=s$ for some $s\geq0$ and set
$p+q+r=s+1$ we then get
\begin{align}
[[A,B],C]\bullet x &= [[A,B],C\bullet x ]+(-1)^{|C||x|}[[A,B]\bullet x,C]\nn\\
&=[A,[B,C\bullet x]]+(-1)^{|B|(|C|+|x|)}[[A,C\bullet x],B]\nn\\
&\quad\,+(-1)^{|C||x|}[[A,B\bullet x],C]+(-1)^{|B||x|+|C||x|}[[A\bullet x,B],C]\nn\\
&=[A,[B,C\bullet x]]+(-1)^{|B|(|C|+|x|)}[[A,C\bullet x],B]\nn\\
&\quad\,+(-1)^{|C||x|}[A,[B\bullet x,C]]+(-1)^{|B||x|+|C||x|}[A\bullet x,[B,C]]\nn\\
&\quad\,+(-1)^{|B||C|}[[A,C],B\bullet x]+(-1)^{|B||C|+|B||x|+|C||x|}[[A\bullet x,C],B]\nn\\
&=[A,[B,C]\bullet x]+(-1)^{|B|(|C|+|x|)}[[A,C]\bullet x,B]\nn\\
&\quad\,+(-1)^{|B||x|+|C||x|}[A\bullet x,[B,C]]\nn\\
&\quad\,+(-1)^{|B||C|}[[A,C],B\bullet x]\nn\\
&=[A,[B,C]]\bullet x+(-1)^{|B||C|}[[A,C],B]\bullet x
\end{align}
using the induction hypothesis in the second and third steps, and the proposition follows by the principle of induction.
\qed

\noindent
It is furthermore easy to see that this Lie algebra is $\mathbb{Z}$-graded, but the $\mathbb{Z}$-grading is different from the one that
is respected by $\dlb \cdot,\cdot \drb$ (on the subspace $U_{0-}$). We have
\begin{align} \label{skiftadgradering}
[\mathcal U_{-p+1},\mathcal U_{-q+1}]\subseteq \mathcal U_{-(p+q)+1}
\end{align}
so the relevant $\mathbb{Z}$-degree of an operator is just (the negative of) its order. 

We will now show that the subspace $\mathcal S$ of the Lie superalgebra $\mathcal U$
closes under the bracket (\ref{bracketdef}) and thus form a subalgebra.
\begin{prop} \label{symmclosure}
If $A,B \in \mathcal S$, then $[A,B] \in \mathcal S$ as well.
\end{prop}
\Pf 
Since all operators of order zero or one are included in $\mathcal S$, and because of the $\mathbb{Z}$-grading (\ref{skiftadgradering}),
we can assume that both $A$ and $B$ are of order one or higher, so that $[A,B]$ is of order two or higher.

We have to show that $\dlb\dlb A,B\drb,u\drb = 0$ for any $u \in U_{2+}$. We first show this for $u \in U_{2+}$, and in particular when $u=\dlb x,y\drb$
for $x,y\in U_1$.
Thus we have to show that
\begin{align}
\big([A,B] \bullet x\big)\bullet y - (-1)^{|x||y|} \big([A,B] \bullet y\big)\bullet x =0 \label{uttryck}
\end{align}
under the assumption that
\begin{align} \label{assA}
(A \bullet x)\bullet y - (-1)^{|x||y|} (A \bullet y)\bullet x =0
\end{align}
and
\begin{align} \label{assB}
(B \bullet x)\bullet y - (-1)^{|x||y|} (B \bullet y)\bullet x = 0 \,.
\end{align}
The first term in (\ref{uttryck}) is equal to 
\begin{align}
[A,B \bullet x] \bullet y +(-1)^{|B||x|} [A \bullet x, B]\bullet y\
&=[A,(B\bullet x)\bullet y]\nn\\&\quad\,+(-1)^{|y|(|B|+|x|)}[A\bullet y,B\bullet x]\nn\\
&\quad\,+(-1)^{|B||x|}[A \bullet x,B\bullet y]\nn\\
&\quad\,(-1)^{|B||x|+|B||y|}[(A\bullet x)\bullet y,B]\,.
\end{align}
Now the second and third term on the right hand side cancel the corresponding contributions from the second term in (\ref{uttryck}). 
Furthermore, the first and fourth term cancel the corresponding contributions from the second term in (\ref{uttryck}) by (\ref{assA}) and (\ref{assB}).
When $u \in U_k$ for $k\geq 3$, we can assume $u=\dlb x, v\drb$, where $v \in U_{2+}$. If we then assume that $\dlb [A,B], v \drb=0$ (as induction hypothesis),
we get
\begin{align}
\dlb [A,B],\dlb x,v\drb \drb &= \dlb \dlb[A,B], x\drb,v \drb %\nn\\
%&= 
=\dlb [A,B]\circ x,v \drb\,,
\end{align}
which is proportional to
\begin{align}
\dlb [A,B]\bullet x,v \drb = \dlb [A,B\bullet x],v \drb +(-1)^{|B||x|} \dlb [A\bullet x,B],v \drb\,.
\end{align}
Now, since $A$ and $B$ are symmetric, $A\bullet x$ and $B\bullet x$ are symmetric as well, and the proposition can be proven by induction.
\qed

\noindent
For $A_p \in \mathcal S_{-p+1}$ and $B_q \in \mathcal S_{-q+1}$, considering the operator $[A_p,B_q]$ as a linear map $(\mathcal U_1)^{p+q}\to \mathcal U_1$,
we have 
\begin{align}
[A_p,B_q](x_1,\ldots,x_{p+q})=[A_p(x_{\langle 1},\ldots,x_{p}),B_q( x_{p+1},\ldots,x_{q \rangle})]\,.
\end{align}

The next proposition says that the identity (\ref{bracketdef}) can be generalised in the sense that $x \in \mathcal U_1$ can be replaced by any
$C \in \mathcal  U_{1-}$. We omit the proof since the steps are the same as in Proposition \ref{u1-isaliealgebra}.
\begin{prop}
For any $A,B,C \in \mathcal  U_{1-}$, we have
\begin{align}
[A,B] \bullet C = [A,B\bullet C]+(-1)^{|B||C|}[A \bullet C,B]\,.
\end{align}
\end{prop}

\subsection{Multiple brackets involving the identity map}

The identity map on $\mathcal U_1$ is an even symmetric operator of order one. 
We denote it simply by $1$, so that $1 \bullet x = 1(x)=x$ and
\begin{align}
[A,1] \bullet x=[A,1 \bullet x]+[A \bullet x,1]=[A,x]+[A \bullet x,1]\,.
\end{align}
We now generalise this notation and, for any integer $k\geq1$ write
\begin{align} \label{111}
[A,k]=[\cdots[[A,1],1],\ldots,1]\,
\end{align}
where the identity map 1 appears $k$ times on the right hand side.
We will furthermore from now on use multibrackets to denote nested brackets (for any elements in any Lie superalgebra) and write
(\ref{111}) as
\begin{align}
[A,k]=[A,1,1\ldots,1]\,.
\end{align}
Note that $[A,i,j]=
[A,i+j]$.

\begin{prop}
If $A \in \mathcal S$, then $[A,k]\in \mathcal S$ as well.
\end{prop}
\Pf
By induction, using that $[A,k,1]=[A,k+1]$, it suffices to show this in the case when $k=1$,
which is a special case of Proposition \ref{symmclosure}.
\qed

\noindent
For $A_p \in \mathcal S_{-p+1}$, considering $[A_p,q]$ as linear map $\mathcal U_1{}^{p+q} \to \mathcal U_1$, we have
\begin{align} \label{vinkelnotation}
[A_p,q](x_1,\ldots,x_{p+q})=[A_p(x_{\langle x_1},\ldots,x_p),x_{p+1},\ldots,x_{p+q \rangle}]\,.
\end{align}

In calculations with multiple brackets involving the identity map, we will need the rules in the next proposition.
They are more or less obvious when reformulated in the notation (\ref{vinkelnotation}) and also straightforward
to prove rigorously
in the more compact notation that we have demonstrated here.
However, since the calculations are rather lengthy, and we have already given similar proofs, we omit this one.

\begin{prop} \label{lemma2prop}
Let $A$ and $B$ be operators and $n\geq 1$ an integer. Then we have
\begin{align}
[A,B,n] = \sum_{k=0}^n\binom{n}{k}[A,k,[B,n-k]]
\end{align}
and
\begin{align} \label{lemma2}
[A,n] \bullet B &= \sum_{i+j=n-1} [A,i,B,j]+[A \bullet B,n]\,\nn\\
&= \sum_{i+j=n-1}\binom{n}{j+1} [A,i,[B,j]]+[A \bullet B,n]\,.
\end{align}
\end{prop}

\noindent
In the summations in (\ref{lemma2}), the summation variables $i$ and $j$ take all non-negative integer values (such that $i+j=n-1$),
and we set $[A,0]=A$ for any operator $A$. Also in all summations
below, the summation variables are allowed to be zero, unless otherwise stated.

%\newpage

\section{Main theorem}

Suppose that the Lie superalgebra $\mathcal U_1=\alg$ 
decomposes into a direct sum $\alg=\ssub\oplus\scom$ of a subalgebra $\ssub$
and a subspace $\scom$.
For any $x \in \alg$, we write $x = x_\sub+x_{\com}$, where $x_\sub \in {\ssub}$ and $x_\com \in \scom$.
Since $\ssub$ is a subalgebra, we thus have
\begin{align}
[a_\sub,b_{\sub}]=[a_\sub,b_\sub]_\sub
\end{align}
and $[a_\sub,b_\sub]_\com=0$ for any $a,b\in \alg$.
Then this decomposition of $\mathcal U_1=\alg$ extends to a decomposition of the Lie superalgebra $\mathcal S$
into a corresponding direct sum $\mathcal S= \mathcal S_\sub \oplus \mathcal S_\com$, where $\mathcal S_\sub$ is a subalgebra. For any $A \in \mathcal U_{1-}$,
we define
$A_\sub$ and $A_\com$
recursively by
\begin{align} \label{Aproj}
A_\sub \bullet x &= (A \bullet x)_\sub\,, & A_\com \bullet x &= (A \bullet x)_\com
\end{align}
for any $x \in \alg$. It follows immediately that $A=A_\sub + A_\com$, and also that  
if $A \in \mathcal S$,
then $A_\sub \in \mathcal S$ and $A_\com \in \mathcal S$ as well.
We let $\mathcal S_\sub$ and $\mathcal S_\com$ be the subspaces spanned by all $A_\sub$ and $A_\com$, respectively,
such that $A \in \mathcal S$. We then have the following proposition, which can be proven in the same way as Proposition \ref{u1-isaliealgebra}, by induction over the
sum of the orders of $A$ and $B$.
\begin{prop}
For any $A,B \in \mathcal S$, we have
\begin{align}
[A_\sub,B_\sub] = [A_\sub,B_\sub]_\sub\,.
\end{align}
Thus the subspace $\mathcal S_\sub$ of the Lie superalgebra $\mathcal S$ is a subalgebra.
\end{prop}
The next proposition says that the identity (\ref{Aproj}) can be generalised in the sense that $x \in \mathcal U_1$ can be replaced by
any $B \in \mathcal S$. Again, it can be proven in the same way as Proposition \ref{u1-isaliealgebra}, by induction over the order of $B$.
\begin{prop} \label{ytterligareenprop}
For any $A,B \in \mathcal S$ we have
\begin{align}
A_\sub \bullet B&= (A \bullet B)_\sub\,, & A_\com \bullet B&= (A \bullet B)_\com
\end{align} 
\end{prop}
We thus obtain the following generalisation of Proposition \ref{lemma2prop} by projecting all outermost brackets on $\scom$.
\begin{prop} \label{lemma2propproj}
Let $A$ and $B$ be operators and $n\geq 1$ an integer. Then we have
\begin{align}
[A,B,n]_\com = \sum_{k=0}^n\binom{n}{k}[A,k,[B,n-k]]_\com
\end{align}
and
\begin{align} \label{lemma2proj}
[A_p,n]_\com \bullet B_q &= \sum_{i+j=n-1} [A_p,i,B_q,j]_\com+[A_p \bullet B_q,n]_\com\,\nn\\
&= \sum_{i+j=n-1}\binom{n}{j+1} [A_p,i,[B_q,j]]_\com+[A_p \bullet B_q,n]_\com\,.
\end{align}
\end{prop}
\Pf This follows directly from Propositions \ref{lemma2prop} and \ref{ytterligareenprop}.
\qed

%\newpage

\noindent
In particular, when $n=1$ we have the identity
\begin{align} \label{litetlemma}
[A,1]_\com \bullet B = [A,B]_\com + [A \bullet B,1]_\com\,,
\end{align}
which we will use below (in the case where $A$ and $B$ are operators of order zero, so that the second term vanishes).

For any $a \in \alg$ and any integer $p\geq0$, we define $a(p) \in \mathcal S_{1-p}$ recursively by
\begin{align}
a(0) &= a_\com\,
\end{align}
and
\begin{align}
a(p) &=\frac1{p!}[a,p]_\com - \sum_{q+r=p-1}\frac1{(r+2)!}[a(q),r+1]_\com
\end{align}
for $p\geq1$.
In particular, we have 
\begin{align} \label{a1}
a(1)=[a,1]_\com-\frac12[a_\com,1]_\com\,,
\end{align}
and thus
\begin{align}
a(1)
\bullet x = [a,x]_\com-\frac12[a_\com,x]_\com\,.
\end{align}
For example,
\begin{align}
a(2)(x_1,x_2)&=\frac1{2!}[a,x_{\langle 1},x_{2\rangle}]_\com-\frac1{2!}[[a,x_{\langle 1}]_\com,x_{2\rangle}]_\com\nn\\
&\quad\,-\frac1{3!}[a_\com,x_{\langle 1},x_{2\rangle}]_\com+\frac1{2!2!}[[a_\com,x_{\langle 1}]_\com,x_{2\rangle}]_\com\,.
\end{align}
For any $a \in \alg$ we also
define 
\begin{align} \label{tildedef}
\tilde a(p) &=\frac1{p!}[a,p] - \sum_{q+r=p}\frac1{(r+1)!}[a(q),r]\,. 
\end{align}
This is however not a recursive definition, since it is $a(q)$, not $\tilde a(q)$, that appears in the second term, and $r+1 \geq 1$ is replaced
by $r\geq0$.
Note also that the bracket is not projected on $\scom$. In fact, $\tilde a(p)$ is projected on $\ssub$, since 
\begin{align}
\tilde a(0) &= [a,0]-[a(0),0]=a-a(0)=a-a_\com = a_\sub
\end{align}
and
\begin{align}
\tilde a(p) &= \frac1{p!}[a,p]
- \sum_{q+s=p-1}\frac1{(s+2)!}[a(q),s+1]-a(p)\nn\\
&=\frac1{p!}[a,p] 
- \sum_{q+s=p-1}\frac1{(s+2)!}[a(q),s+1]\nn\\
&\quad\, -\frac1{p!}[a,p]_\com + \sum_{q+r=p-1}\frac1{(r+2)!}[a(q),r+1]_\com \nn\\
&= \frac1{p!}[a,p]_\sub - \sum_{q+r=p-1}\frac1{(r+2)!}[a(q),r+1]_\sub
\end{align}
for $p\geq 1$.
It follows that
\begin{align}
[\tilde a(p),\tilde b(q)]_\com = 0\, \label{tildevanish}
\end{align}
for any $a,b$ and $p,q \geq0$.

We are now ready to formulate and prove our main theorem.

\begin{theorem}\label{thetheorem}
The map
\begin{align}
\alg \to \mathcal{S}_\com\,, \quad a \mapsto \sum_{p=0}^\infty a(p)
\end{align}
is a Lie superalgebra homomorphism.
\end{theorem}
\Pf We will show that
\begin{align}
\blb \sum_{p=0}^\infty a(p), \sum_{q=0}^\infty b(q) \brb = \sum_{r=0}^\infty [a,b](r)
\end{align}
for any $a,b \in \alg$. The left hand side is equal to
\begin{align}
\sum_{p=0}^\infty\sum_{q=0}^\infty \blb  a(p), b(q) \brb = \sum_{r=0}^\infty\ \sum_{{p+q=r+1}} \blb  a(p), b(q) \brb\,.
\end{align}
Thus it suffices to show that
\begin{align}
\sum_{{p+q=r+1}} \blb  a(p), b(q) \brb = [a,b](r) \label{indlikhet}
\end{align}
for $r=0,1,2,\ldots$. We will do this by induction. When $r=0$, the left hand side in (\ref{indlikhet}) equals
\begin{align}
\sum_{{p+q=1}} \blb  a(p), b(q) \brb &= \blb  a(0), b(1) \brb+\blb  a(1), b(0) \brb\nn\\
&=\blb a_\com,[b,1]_\com\brb - \frac12\blb a_\com,[b_\com,1]_\com\brb\nn\\
&\quad\,+\blb[a,1]_\com,b_\com\brb-\frac12\blb[a_\com,1]_\com,b_\com\brb\,\nn\\
&=-(-1)^{|a||b|}[b,1]\bullet a_\com+\frac12(-1)^{|a||b|}[b_\com,1]_\com \bullet a_\com\nn\\
&\quad\,+[a,1]_\com\bullet b_\com-\frac12[a_\com,1]_\com\bullet b_\com\nn\\
&=-(-1)^{|a||b|}[b,a_\com]+\frac12(-1)^{|a||b|}[b_\com,a_\com]_\com\nn\\
&\quad\,+[a,b_\com]_\com-\frac12[a_\com,b_\com]_\com\nn\\
&=[a,b_\com]_\com+[a_\com,b]_\com-[a_\com,b_\com]_\com\,,
\end{align}
where we have used (\ref{litetlemma}) and (\ref{a1}),
whereas the right hand side
equals $[a,b]_\com$. Thus the right hand side minus the left hand side equals
\begin{align}
[a,b]_\com-[a,b_\com]_\com-[a_\com,b]_\com+[a_\com,b_\com]_\com=[a-a_\com,b-b_\com]_\com=[a_\sub,b_\sub]_\com=0\,.
\end{align}

In the induction step, we need to study 
\begin{align} \label{sum}
\sum_{m+n=k} \blb a(m),b(n) \brb  &= 
\sum_{m+n=k} \bigg(a(m) \bullet b(n) -(-1)^{|a||b|} b(n) \bullet a(m)\bigg)
\end{align}
for some $k\geq2$
and show that this the expression equals
\begin{align}
[a,b](k-1)
\end{align}
under the assumption that
\begin{align}
\sum_{m+n=s} \blb a(m),b(n) \brb  &= [a,b](s-1)
\end{align}
for $s=1,\ldots,k-1$.

We will first study the first term in the summand on the right hand side of (\ref{sum}), and a particular part of it. Its counterpart, the corresponding part of the second term in the summand, is then obtained by interchanging $a$ and $b$, and multiplying with $(-1)^{|a||b|}$.
In the summations in (\ref{517}) and (\ref{518}) below where the summation variables add up to $m-1$, the sum should be read as zero if $m=0$.

%\newpage
\noindent
We have
\begin{align}
a(m) \bullet b(n)&=
\left(\frac1{m!}[a,m]_\com-\sum_{p+q=m-1}\frac1{(q+2)!}[a(p),q+1]_\com\right)\bullet b(n)\nn\\
&=
\frac1{m!}[a,m]_\com\bullet b(n)-\sum_{p+q=m-1}\frac1{(q+2)!}[a(p),q+1]_\com\bullet b(n)\nn\\
&=
\frac1{m!}[a \bullet b(n) ,m]_\com +\frac1{m!}\sum_{p+q=m-1}\binom{m}{q+1}\big[[a,p],[b(n),q]\big]_\com   \nn\\
&\quad\,-\sum_{p+q=m-1}\frac1{(q+2)!}[a(p) \bullet b(n),q+1]_\com\nn\\
&\quad\,-\sum_{p+q=m-1}\frac1{(q+2)!}\sum_{r+s=q}\binom{q+1}{s+1}\big[[a(p),r],[ b(n),s]\big]_\com\nn\\
&=
\sum_{p+q=m-1}\frac1{p!(q+1)!}\big[[a,p],[b(n),q]\big]_\com   \nn\\
&\quad\,-\sum_{p+q=m-1}\frac1{(q+2)!}[a(p) \bullet b(n),q+1]_\com\nn\\
&\quad\,-\sum_{p+r+s=m-1}\frac1{(r+s+2)!}\binom{r+s+1}{s+1}\big[[a(p),r],[b(n),s]\big]_\com\,. \label{517}
\end{align}
The contribution to the sum (\ref{sum}) from the last term in (\ref{517}), and its counterpart, is
\begin{align}
\sum_{m+n=k}&
\Bigg(-\sum_{p+r+s=m-1}\frac1{(r+s+2)!}\binom{r+s+1}{s+1}\big[[a(p),r],[b(n),s]\big]_\com\nn\\
&\quad+(-1)^{|a||b|} \sum_{q+r+s=n-1}\frac1{(r+s+2)!}\binom{r+s+1}{r+1}\big[[b(q),s],[a(m),r]\big]_\com
\Bigg)\nn\\
&=\sum_{p+q+r+s=k-1}\Bigg(-\frac1{(r+s+2)!}\binom{r+s+1}{s+1}\big[[a(p),r],[b(q),s]\big]_\com\nn\\
&\qquad\qquad\qquad\qquad\!\!
+(-1)^{|a||b|} \frac1{(r+s+2)!}\binom{r+s+1}{r+1}\big[[b(q),s],[a(p),r]\big]_\com\Bigg)\nn\\
&=-\sum_{p+q+r+s=k-1}\frac1{(r+s+2)!}\Bigg(\binom{r+s+1}{s+1}+\binom{r+s+1}{r+1}\Bigg)\big[[a(p),r],[b(q),s]\big]_\com\nn\\
&=- \sum_{p+q+r+s=k-1}\frac1{(r+s+2)!}\binom{r+s+2}{r+1}\big[[a(p),r],[b(q),s]\big]_\com\nn\\
&=- \sum_{p+q+r+s=k-1}\frac1{(r+1)!}\frac1{(s+1)!}\big[[a(p),r],[b(q),s]\big]_\com\,. \label{518}
\end{align}
Taking all terms into account, we get
\begin{align}
\sum_{m+n=k} \blb a(m),b(n) \brb  &= 
\sum_{m+n=k} \bigg(a(m) \bullet b(n) -(-1)^{ab} b(n) \bullet a(m)\bigg)\nn\\
&=
\sum_{q+r+s=k-1}\frac1{r!(s+1)!}\big[[a,r],[b(q),s]\big]_\com  \tag{a} \\
&\quad\,-\sum_{p+q+r=k-1}\frac1{(r+2)!}[a(p) \bullet b(q),r+1]_\com\nn\\
&\quad\,-\sum_{p+q+r+s=k-1}\frac1{(r+s+2)!}\binom{r+s+1}{s+1}\big[[a(p),r],[b(q),s]\big]_\com\nn\\
&\quad\,-(-1)^{ab}\sum_{p+r+s=k-1}\frac1{s!(r+1)!}\big[[b,s],[a(p),r]\big]_\com   \nn\\
&\quad\,+(-1)^{ab}\sum_{p+q+r=k-1}\frac1{(r+2)!}[b(q) \bullet a(p),r+1]_\com\nn\\
&\quad\,+(-1)^{ab}\sum_{p+q+r+s=k-1}\frac1{(r+s+2)!}\binom{r+s+1}{r+1}\big[[b(q),s],[a(p),r]\big]_\com\nn\\
&= \sum_{q+r+s=k-1}\frac1{r!(s+1)!}\big[[a,r],[b(q),s]\big]_\com \tag{b}\\
&\quad\,+\sum_{p+r+s=k-1}\frac1{s!(r+1)!}\big[[a(p),r],[b,s]\big]_\com\nn\\
&\quad\,- \sum_{p+q+r=k-1}\frac1{(r+2)!}\big[\blb a(p) , b(q) \brb,r+1\big]_\com\nn\\
&\quad\,- \sum_{p+q+r+s=k-1}\frac1{(r+1)!}\frac1{(s+1)!}\big[[a(p),r],[b(q),s]\big]_\com\,\nn\\
&= -\sum_{p+q=k-1}[\tilde a(p),\tilde b(q)]_\com+\sum_{p+q=k-1}\frac1{p!\,q!}\big[[a,p],[b,q]\big]_\com \tag{c}\\
&\quad\,- \sum_{p+q=k-2}\frac1{(q+2)!}[[a,b](p) ,q+1]_\com\nn\\
&= \frac{1}{(k-1)!}\sum_{p+q=k-1}\binom{k-1}{q}\big[[a,p],[b,q]\big]_\com \tag{d}\\
&\quad\,
- \sum_{p+q=k-2}\frac1{(q+2)!}[[a,b](p) ,q+1]_\com\nn\\
&= \frac{1}{(k-1)!}\big[[a,b],k-1\big]_\com - \sum_{p+q=k-2}\frac1{(q+2)!}[[a,b](p) ,q+1]_\com\tag{e}\\
&=[a,b](k-1)\,.
\end{align}
Here we have used (\ref{517}) in (a). In (b) we have used the definition of $\blb\cdot,\cdot\brb$
and (\ref{518}). The second and the fifth term on the right hand side of (a) go into the
third term of the right hand side of (b), whereas the third and sixth term of 
(a) go into the
fourth term of the right hand side of (b), by (\ref{518}). In (c) we have used the definition (\ref{tildedef}) of
$\tilde a (p)$ and $\tilde b (q)$, and the induction hypothesis. In (d) we have used (\ref{tildevanish})
and in (e) we have used Proposition \ref{lemma2propproj}.
The theorem now follows by the principle of induction.
\qed

\noindent
Considering $a(p)=a_p$ as a linear map $\alg^p \to \scom$, we have  
\begin{align}
a_p(x_1,\ldots,x_p) &= \sum_{k=0}^p \sum \frac1{m_1!}\frac{(-1)}{(m_2-m_1+1)!}\cdots \frac{(-1)}{(p-m_{k}+1)!}\times\nn\\
&\quad\quad\quad\quad\times [[\cdots [[a,x_{\langle 1},\ldots,x_{m_1}]_\com,
x_{m_1+1}, \ldots, x_{m_2}]_\com ,x_{m_2+1},\ldots\nn\\
&\qquad\qquad\qquad\qquad\qquad\ \,\ldots, 
x_{m_k}]_\com, x_{m_k+1} \ldots, x_{p\rangle}]_\com\,,\label{genform} 
\end{align}
where the
inner sum goes over all $k$-tuples of integers $(m_1,\ldots,m_k)$ such that
\begin{align}
0 \leq m_1 < m_2 < \cdots < m_k < p\,.
\end{align}
If $m_1=0$ (and $k>0$), the factor in the second and third line should be read as
\begin{align}
&\quad\quad\times [[\cdots [a_\com ,
x_{\langle 1}, \ldots, x_{m_2}]_\com ,x_{m_2+1},\ldots\nn\\
&\qquad\qquad\qquad\qquad\ \, \, \ldots, 
x_{m_k}]_\com, x_{m_k+1} \ldots, x_{p\rangle}]_\com\,.
\end{align}
If $k=0$, the inner sum in (\ref{genform}) should be read as
\begin{align}
\frac1{p!}[a,x_{\langle 1}, \ldots, x_{p\rangle}]_\com\,.
\end{align}
Here $x_1,\ldots,x_p \in \alg$, but since $\scom$ is a subspace of $\alg$, we can as well assume $x_1,\ldots,x_p \in \scom$
and consider $a(p)=a_p$ as a linear map $\scom^p \to \scom$.

\section{Getzler's theorem}

As an example of an application, we end this paper by proving a theorem which says that any
differential graded
Lie algebra (a Lie superalgebra with a consistent $\mathbb{Z}$-grading and a differential) gives rise to an $L_\infty$-algebra (a generalisation of a
differential graded
Lie algebra including also higher brackets \cite{Lada:1992wc,Marklada}).
Combined with the result described in Section 3.1 that any Leibniz algebra gives rise to a differential graded
Lie algebra, it leads to the conclusion that any Leibniz algebra gives rise to an $L_\infty$-algebra \cite{Kotov:2018vcz,Lavau:2019oja,Lavau:2020pwa}.
The theorem has already been proven in at least two different ways in the literature. It follows from the results in Ref.\ \cite{FiorenzaManetti} by Fiorenza and Manetti,
and has been proven more directly in Ref.\ \cite{Getzler:2010} by Getzler. Here we follow Getzler's formulation of the it,
and prove it using our main result, Theorem \ref{thetheorem}.

Suppose that the Lie superalgebra $\alg$ has a consistent $\mathbb{Z}$-grading,
$\alg = \bigoplus_{i\in\mathbb{Z}}\alg^{(i)}$.
Then this
$\mathbb{Z}$-grading induces a $\mathbb{Z}$-grading on each subspace $\mathcal S_i$ of $\mathcal S$, and thus a $\mathbb{Z}$-grading
of $\mathcal S$, different from the one
that $\mathcal S$ comes with by construction,
\begin{align}
\mathcal S_i &= \bigoplus_{j\in\mathbb{Z}}\mathcal S_i{}^{(j)}\,,&
\mathcal S^{(j)} &= \bigoplus_{i\in\mathbb{Z}}\mathcal S_i{}^{(j)}\,.
\end{align}
The two $\mathbb{Z}$-gradings form together a $(\mathbb{Z}\times \mathbb{Z})$-grading,
\begin{align}
\mathcal S = \bigoplus_{(i,j)\in\mathbb{Z}\times \mathbb{Z}}\mathcal S_i{}^{(j)}= \bigoplus_{i\in\mathbb{Z}}\mathcal S_i =
\bigoplus_{j\in\mathbb{Z}}\mathcal S^{(j)}\,.
\end{align}
If there is an element $Q\in \mathcal S_{0-}{}^{(-1)}$ such that $\blb Q,Q\brb=0$, then $\alg$ together with $Q$ constitutes an $L_\infty$-{algebra}.
The element $Q$ can then be decomposed
as a sum of elements $Q_{p} \in S_{-p+1}$, for $p=1,2,\ldots$,
each of which
can be considered as a linear map $\alg^p \to \alg$, called a $p$-{\it bracket}. Following Ref.\ \cite{Getzler:2010}, we use 
curly brackets for the $p$-brackets, $Q_{p}(x_1,\ldots,x_p)=\{x_1,\ldots,x_p\}$.
The condition $\blb Q,Q\brb=0$ decomposes into infinitely many identities for these $p$-brackets, similar to the usual Jacobi identity.

We note that there are different conventions for $L_\infty$-algebras. The fact that we consider the $p$-brackets as elements $Q_p\in \mathcal S_{-p+1}{}^{(-1)}$ means that
we use the convention where they are graded symmetric
and have degree $-1$.

\begin{theorem}{\rm\cite{Getzler:2010}}
Let $L=\bigoplus_{i \in \mathbb{Z}} L_i$ be a differential graded Lie algebra with differential $\delta$ of degree $-1$ and bracket $[\cdot,\cdot]$. Let $D$
be the linear operator on $L$ which equals $\delta$ on $L_1$ but vanishes on $L_i$ for $i\neq 1$. Then the subspace $\bigoplus_{i \geq 1}L_i$ is
an $L_\infty$-algebra with $p$-brackets given by
\begin{align} \label{1-bracket}
\{ x \} = \delta x - Dx 
\end{align}
for $p=1$ and
\begin{align} \label{2-bracketsochhogre}
\{x_1,x_2,\ldots,x_p\}=
-\frac1{(p-1)!}B_{p-1}{}^-
[Dx_{\langle1},x_{2},\ldots,x_{p\rangle}]
\end{align}
for $p=2,3,\ldots$ where $B_n{}^-$ are the Bernoulli numbers ($B_n{}^-=-\frac12, \frac16,0,-\frac1{30},0,\ldots$ for $n=1,2,3,4,5,\ldots)$.
\end{theorem}

\noindent
In Ref.~\cite{Getzler:2010} there appears to be a sign error that we have here corrected by inserting a minus sign on the right hand side
of (\ref{2-bracketsochhogre}) \cite{Lavau:2019oja}. The occurrence of Bernoulli numbers in this context was observed in Ref.~\cite{Bering:2006eb}, and 
it was also shown in Ref.~\cite{Cederwall:2018aab} that they similarly show up in extended geometry,
encoding the gauge structure of generalised diffeomorphisms.\\

\Pf Let $\mathcal G_{-1}$ be a one-dimensional vector space spanned by an element $\Theta$ and set $\mathcal G_i=L_i$ for $i=0,1,2,\ldots$.
Then 
\begin{align}
\mathcal G=\mathcal G_{-1}\oplus \mathcal G_{0} \oplus \mathcal G_{1} \oplus \cdots
\end{align}
is a consistently $\mathbb{Z}$-graded Lie superalgebra, 
where the bracket in the subalgebra $\bigoplus_{i\geq 0} L_i$ of $L$ is extended by $[\Theta,\Theta]=0$ and $[\Theta,x]=\delta x$ for
$x \in \mathcal G_{0},\mathcal G_{1},\ldots$.
Furthermore, set
$\ssub=\mathcal G_{-1}\oplus \mathcal G_{0}$ and $\scom = \mathcal G_{1}\oplus \mathcal G_{2}\oplus \cdots$.
Then $\ssub$ is a subalgebra of $\mathcal G$ and
$\mathcal G=\ssub\oplus \scom$.
We can thus use Theorem \ref{thetheorem}, which
in particular says that the element $Q=\sum_{p=0}^\infty \Theta(p)$ in $\mathcal S$
satisfies $\blb Q,Q \brb=0$, since $[\Theta,\Theta]=0$ in $\alg$.
Also, it follows by the construction of $\Theta(p)$ and the 
$\mathbb{Z}$-grading that $\Theta(p) \in \mathcal S^{(-1)}$, since $\Theta \in \mathcal S^{(-1)}$ and the identity map
$1 \in \mathcal S^{(0)}$. 
If we write $\Theta(p)=Q_p$,
it thus follows that $\scom$ together with the element $Q=\sum_{p=1}^\infty Q_p$ in $\mathcal S$ (note that $Q_0=\Theta(0)=\Theta_\com=0$)
is an $L_\infty$-algebra,
with the $p$-brackets
\begin{align}
\{x_1,\ldots,x_p\}=Q_p(x_1,\ldots,x_p)\,.
\end{align}
The right hand side here is 
given by (\ref{genform}) with $a=\Theta$, that is
\begin{align}
Q_p(x_1,\ldots,x_p) &= \sum_{k=0}^p \sum \frac1{m_1!}\frac{(-1)}{(m_2-m_1+1)!}\cdots \frac{(-1)}{(p-m_{k}+1)!}\times\nn\\
&\quad\quad\quad\quad\times [[\cdots [[\Theta,x_{\langle 1},\ldots,x_{m_1}]_\com ,
x_{m_1+1}, \ldots, x_{m_2}]_\com ,x_{m_2+1},\ldots\nn\\
&\qquad\qquad\qquad\qquad\qquad\ \,\,\ldots, 
x_{m_k}]_\com, x_{m_k+1} \ldots, x_{p\rangle}]_\com\,,\label{genformtheta} 
\end{align}
where the inner sum goes over all $k$-tuples of integers $(m_1,\ldots,m_k)$ such that
\begin{align}
0 \leq m_1 < m_2 < \cdots < m_k < p\,.
\end{align}

It remains to show that (\ref{genformtheta}) equals the expressions on the right hand side of (\ref{1-bracket}) and (\ref{2-bracketsochhogre}) when $p=1$ and $p\geq 2$, respectively.
When $p=1$,
we indeed get
\begin{align}
\Theta (x) = [\Theta,x]_\com -\frac12 [\Theta_\com,x]_\com = [\Theta,x]_\com = (\delta x)_\com = \delta x-Dx\,.
\end{align}
When $p\geq 2$, all terms in (\ref{genformtheta}) with $m_1=0$ are zero, since $\Theta_\com=0$. Furthermore, all the subscripts ${}_\com$ but the first one can be removed, since
$\scom$ is a subalgebra in this case. Also, when $m_1\geq 2$ even the first subscript ${}_\com$ can be removed.
Thus, for $p\geq 2$ we have
\begin{align}
Q_p(x_1,\ldots,x_p) &= \sum_{k=0}^p \sum \frac1{1!}\frac{(-1)}{n_1!}\frac{(-1)}{(n_2-n_1+1)!}\cdots \frac{(-1)}{(p-n_{k}+1)!}\times\nn\\
&\quad\quad\quad\quad\times [[\cdots [[\Theta,x_{\langle 1}]_\com,
x_{2}, \ldots, x_{n_1}],x_{n_1+1}, \ldots, x_{n_2}] ,x_{n_2+1},\ldots\nn\\
&\qquad\qquad\qquad\qquad\qquad\qquad \,\,\,\ldots, 
x_{n_k}], x_{n_k+1} \ldots, x_{p\rangle}]\,\nn\\
&\quad+\sum_{k=0}^p 
\sum 
\frac1{n_1!}\frac{(-1)}{(n_2-n_1+1)!}\cdots \frac{(-1)}{(p-n_{k}+1)!}\times\nn\\
&\quad\quad\quad\quad\quad\times [[\cdots [[\Theta,x_{\langle 1},\ldots,x_{n_1}],x_{n_1+1}, \ldots, x_{n_2}] ,x_{n_2+1},\ldots\nn\\
&\qquad\qquad\qquad\qquad\qquad\quad \,\,\,\,\ldots, 
x_{n_k}], x_{n_k+1} \ldots, x_{p\rangle}]
\end{align} 
where
the inner sums go over all $k$-tuples of integers
$(n_1,\ldots,n_k)$ such that
\begin{align}
1 < n_1 < n_2 < \cdots < n_k < p\,.
\end{align}
Since we have removed the subscripts ${}_\com$, the factor in the second and third line of each term above is actually independent of the choice of 
$n_1,n_2,\ldots,n_k$, and also of the integer $k$.
Setting
\begin{align}
C_p=\sum_{k=0}^{p}
\sum  \frac{1}{n_1!}\frac{(-1)}{(n_2-n_1+1)!}\cdots \frac{(-1)}{(p-n_{k}+1)!}\,,
\end{align}
where the term with $k=0$ should be read as $1/p!$,
we thus get
\begin{align}
Q_p(x_1,\ldots,x_p) &= -C_p[[\Theta,x_{\langle 1}]_\com, x_2,\ldots,x_{p\rangle}]\nn\\
&\quad\, + 
 C_p [\Theta,x_{\langle 1},x_{2},\ldots,x_{p\rangle}]\nn\\
 &=-C_p\big([\delta x_{\langle 1},x_{2},\ldots,x_{p\rangle}]-[Dx_{\langle 1},
 x_{2}, \ldots, x_{p\rangle}]\big)\nn\\
 &\quad\,+C_p[\delta x_{\langle 1},x_{2},\ldots,x_{p\rangle}]\nn\\
 &=C_p
  [Dx_{\langle 1},
x_{2}, \ldots, x_{p\rangle}] \,.
\end{align}
Now $C_p$ can also be written 
\begin{align}
C_p = -\sum_{j=1}^p \sum
\frac{(-1)}{(m_1+1)!} \cdots \frac{(-1)}{(m_j+1)!}
\end{align}
where the inner sum goes over all 
$k$-tuples of positive integers $(m_1,\ldots,m_j)$ such that $m_1+\cdots+m_j=p-1$.
Written this way, it is easily shown 
(by induction, using recursion formulas for the Bernoulli numbers) that
\begin{align}
C_p = -\frac1{(p-1)!}B_{p-1}{}^-
\end{align}
and we arrive at (\ref{2-bracketsochhogre}).
\qed

\begingroup\raggedright

\end{document}